\documentclass[11pt,a4paper]{article}
\usepackage{amsmath}
\usepackage{amssymb}
\usepackage{amsthm}
\usepackage{dsfont}
\usepackage{bm}
\usepackage{url}
\usepackage{appendix}
\usepackage{comment}
\usepackage{authblk}
\usepackage[version=3]{mhchem} 

\usepackage{tikz}
\usepackage{amsfonts} 
\usepackage{color}
\usepackage{booktabs}

\def\RR{{\mathbb R}}

\def\R{{\mathbb R}}

\def\tr{{\rm Tr}}

\def\cM{{\mathcal M}}
\def\cN{{\mathcal N}}

\def\cR{{\mathcal R}}

\def\cT{{\mathcal T}}

\def\cW{{\mathcal W}}

\def\br{{\mathbf r}}

\author[1]{Laurent Vidal}
\affil[1]{CERMICS, Ecole des Ponts and Inria Paris, 6 \& 8 avenue Blaise Pascal, 77455~Marne-la-Vall\'ee, France}
\author[2]{Tommaso Nottoli}
\affil[2]{Dipartimento di Chimica e Chimica Industriale, Universit\`a di Pisa, $\mbox{Via G. Moruzzi 13,}$ 56124 Pisa, Italy}
\author[2]{Filippo Lipparini}
\author[1]{Eric Canc\`es}
\title{Geometric Optimization of Restricted-Open and Complete Active Space Self-Consistent Field Wavefunctions}

\begin{document}

\maketitle

\begin{abstract}
We explore Riemannian optimization methods for Restricted-Open-shell Hartree-Fock (ROHF) and Complete Active Space Self-Consistent Field (CASSCF) methods. After showing that ROHF and CASSCF can be reformulated as optimization problems on so-called flag manifolds, we review Riemannian optimization basics and their application to these specific problems. We compare these methods to traditional ones and find robust convergence properties without fine-tuning of numerical parameters. Our study suggests Riemannian optimization as a valuable addition to orbital optimization for ROHF and CASSCF, warranting further investigation.
\end{abstract}

\section{Introduction}
Orbital optimization is one of the most common task performed in quantum chemistry calculations. It is the numerical problem associated with Hartree-Fock (HF) \cite{Hartree1957} and Kohn-Sham Density Functional Theory (KS DFT) \cite{Kohn65}, as well as a component of Complete Active Space Self-Consistent Field (CASSCF) calculations \cite{Werner1987,Shepard1987,Roos1987}, and is further encountered in orbital optimized post-Hartree Fock methods \cite{scuseria1987optimization,sherrill1998energies}.
The various algorithms that have been proposed to tackle this problem can be grouped into two families: fixed point methods, such as Roothaan's SCF algorithm \cite{Roothan51,Roothaan1960}, and direct optimization methods, such as quadratically convergent optimization strategies \cite{Bacskay1981,Bacskay1982}.
For HF and DFT, the former family is the most commonly employed, due to the existence of very robust implementation that exploit convergence acceleration techniques such as Pulay's Direct Inversion in the Iterative Subspace \cite{Pulay1980,Pulay1982,Hamilton1985} (DIIS), constraint relaxation methods such as the Optimal Damping Algorithm \cite{CancesLeBris2000,Cances2000,Cances2001} (ODA), or more sophisticated related techniques such as E-DIIS \cite{Kudin2002} or A-DIIS \cite{HuYang2010}.
Nevertheless, direct optimization techniques have received quite some attention due to their robustness and due to the possibility of implementing them avoiding dense linear algebra operations (e.g., diagonalization of the Fock matrix). 

Direct minimization techniques for Restricted Open-Shell Hartree-Fock (ROHF) calculations are relatively scarce compared to Self-Consistent Field (SCF) methods, which predominantly feature in quantum chemistry software. Among direct minimization approaches, noteworthy methods include Geometric Direct Minimization (GDM) techniques \cite{Dunietz2002} and the QC-SCF algorithm \cite{nottoli_black-box_2021}. Additionally, the Second-Order SCF (SOSCF) algorithm \cite{chaban1997approximate,neese2000approximate} and the DIIS-GDM method \cite{Dunietz2002,Voorhis2002}, which amalgamate aspects of both SCF and direct minimization strategies, merit mention. Moreover, the CUHF method, as introduced by Tsuchimochi and Scuseria \cite{tsuchimochi2010communication}, can be adapted for ROHF computations through the utilization of a direct minimization procedure designed for Unrestricted Hartree-Fock (UHF) calculations.

A difficulty associated with the formulation and implementation of direct minimization techniques is due to the fact that the quantity that needs to be optimized, such as the molecular orbitals (MO) coefficients, or the density matrix, needs to satisfy nonlinear constraints. In other words, the minimization set is not a vector space, but rather a differentiable manifold.

It is well-known that after discretization in a finite basis set, HF and KS models can be formulated as optimization problems on Stiefel (molecular orbital formalism) or Grassmann (density matrix formalism) manifolds~\cite{Edelman1998}. These formulations lead to enlightening geometric interpretations of the Hartree-Fock and Kohn-Sham equations, and to the design and convergence analysis \cite{CancesLeBris2000a,Levitt2012,Cances2021} of robust and efficient direct minimization algorithms. 
The purpose of this article is to show that Restricted Open-Shell Hartree-Fock (ROHF) and Complete Active Space Self-Consistent Field (CASSCF) methods can be reformulated as optimization problems on so-called flag manifolds \cite{Ye2022}. This allows one to shed new light on the ROHF and CASSCF equations, and the direct minimization algorithms used to solve these problems. 
While the work presented in this article does not lead to game-changing improvements in orbital optimization, we hope that it will provide the community with a set of rigorous tools that can be used for further developments, as the ones recently proposed by some of us for the extrapolation of the SCF density matrix in the context of ab-initio molecular dynamics simulations \cite{polack_approximation_2020,Polack_JCTC_GExt,Pes2023}.

This article is organized as follows. In Section~\ref{sec:theory}, we briefly recall the high-spin ROHF and CASSCF orbital optimization problem in terms 
of both density-matrix (DM) and molecular orbitals (MO) formulations and provide a simple geometric interpretation of the ROHF and CASSCF equations. 
In Section~\ref{sec:optimization_riemmannian}, we review the basic concepts of geometric optimization (Riemannian gradient and Hessian, vector transport, affine connection, geodesic, retraction). In Section~\ref{sec:optimization_flag}, we  discuss more specifically geometric optimization for ROHF and CASSCF, providing also the tools to translate any algorithm formulated in the MO formalism into the DM formalism and viceversa. We then provide geometric interpretations of existing direct minimization algorithms for ROHF and CASSCF, and propose new ones. We also introduce a direct optimization method circumventing the use of virtual orbitals, which is useful for ROHF in large basis sets (i.e. planewaves, finite elements, or wavelets). Numerical results are reported in Section~\ref{sec:numerical_results}.

\section{ROHF and CASSCF\label{sec:theory}}
In this section, we briefly recapitulate the orbital optimization problem for ROHF and CASSCF and introduce the manifolds associated with the MO and DM formalisms. ROHF and CASSCF methods indeed share common features. They both involve 
\begin{itemize}
\item a set of $N_I$ doubly-occupied molecular orbitals, often called internal orbitals;
\item a set of $N_A$ partially-occupied molecular orbitals, often called active orbitals,
\end{itemize}
the latter being orthogonal to the former. Consider a molecular system with $N$ electrons discretized in a basis set of size $\cN_{\rm b}$. We denote the electronic Hamiltonian by 
\begin{equation}\label{eq:electronic_Hamiltonian}
\widehat H_N = - \frac{1}{2} \sum_{i=1}^N \nabla_{\br_i}^2 + \sum_{i=1}^N V_{\rm nuc}(\br_i) + \sum_{1 \le i < j \le N} \frac{1}{|\br_i-\br_j|}.
\end{equation}

\paragraph{Notation} To avoid possible misunderstandings, we clarify here the notation that will be adopted throughout the paper. Let us assume that we have discretized the problem using an orthonormal set of atomic orbitals (AO) $\{\chi_\mu\}_{\mu=1}^{\cN_{\rm b}}$, obtained, for instance, by L{\"o}wdin orthogonalization of a usual Gaussian-type or Slater-type basis. This implies, for the overlap matrix:
\[
S_{\mu\nu} = \langle \chi_\mu | \chi_\nu\rangle = \delta_{\mu\nu}.
\] 
We use greek letters $\mu,\nu,\ldots$ to label AOs.
Molecular orbitals $\{\phi_p\}_{p=1}^{\cN_{\rm b}}$ are written as linear combinations of atomic orbitals 
\[
 \phi_p = \sum_{\mu=1}^{\cN_{\rm b}} C_{\mu p} \chi_\mu
\]
where the coefficient matrix $C \in \R^{\cN_{\rm b}\times \cN_{\rm b}}$ is, due to our choice of orthogonal AOs, an orthogonal matrix, i.e., $CC^T = C^TC = I_{\cN_{\rm b}}$.
We divide the molecular orbitals into three sets: $N_I$ \emph{internal} orbitals, that are always doubly occupied, $N_A$ \emph{active} orbitals, that are singly occupied in ROHF and have varying occupation for CASSCF, and $N_E$ \emph{external} orbitals, that are always empty. We use $i,j,\ldots$ to label internal orbitals, $u,v,\ldots$ to label active orbitals, $a, b, \ldots$ to label external orbitals, and $p, q,\ldots$ for generic ones. 
We call $\Pi^I$ and $\Pi^A$ the orthogonal projectors on the space spanned by internal and active orbitals, respectively, in the AO basis, i.e.,
\begin{equation}\label{eq:projectors}
\Pi^I_{\mu\nu} = \sum_{i=1}^{N_I} C_{
\mu i}C_{\nu i}, \quad \Pi^A_{\mu\nu} = \sum_{u=N_I+1}^{N_I+N_A} C_{
\mu u}C_{\nu u}
\end{equation}
$P \in \R^{\cN_{\rm b}\times \cN_{\rm b}}$ denotes the one-body reduced density matrix (1-RDM) in the AO basis, while $\gamma \in \R^{\cN_{\rm b}\times \cN_{\rm b}}$ denotes the one-body reduced density matrix in the MO basis. We call $m_\alpha$ and $m_\beta$ the number of $\alpha$ and $\beta$ \emph{active} electrons, respectively, so that the total number of electrons is given by $N=2N_I + m_\alpha + m_\beta$. In high-spin ROHF, $m_\alpha = N_A$ and $m_\beta = 0$.
Using these conventions, the high-spin ROHF density matrix is given by
\begin{equation}
\label{eq:PROHF}
P^{\rm ROHF} = 2\Pi^I + \Pi^A    
\end{equation}
where we note that the $\alpha$ and $\beta$ one-body spin density matrices (1-SDM) are given by
\begin{equation}
    \label{eq:ROHFSpinRDM}
    P^{\rm ROHF, \alpha} = \Pi^I + \Pi^A, \quad P^{\rm ROHF, \beta} = \Pi^I
\end{equation}
For CASSCF, the MO 1-RDM $\gamma^{\Psi^{\rm CAS}}$ has a block structure, with
\begin{equation}
\gamma^{\Psi^{\rm CAS}}_{ij} = 2\delta_{ij},\quad \gamma^{\Psi^{\rm CAS}}_{uv} = \langle \Psi^{\rm CAS} | \hat{E}_{uv}|\Psi^{\rm CAS}\rangle,    \quad \gamma^{\Psi^{\rm CAS}}_{ab} = 0,
\end{equation}
where
\[
\hat{E}_{pq} = \hat{a}^\dagger_{p\alpha}\hat{a}_{q\alpha} +
\hat{a}^\dagger_{p\beta}\hat{a}_{q\beta}
\]
is the spin-traced singlet excitation operator, and all other blocks vanishing. 
In the AO basis, $P^{\Psi^{\rm CAS}} = C\gamma^{\Psi^{\rm CAS}} C^T$ satisfies thus the following relation:
\begin{equation}
    2\Pi^I \leq P^{\Psi^{\rm CAS}} \leq 2\Pi^I + \Pi^A.
\end{equation}
To define 1-SDM for CASSCF, we need to introduce the active spin densities $\gamma_{uv}^{\Psi^{\rm CAS},\sigma}$, for $\sigma = \alpha, \beta$ which are defined as follows:
\begin{equation}
\label{eq:CasSpinRDM}
\gamma^{\Psi^{\rm CAS},\sigma}_{ij} = \delta_{ij}, \quad
\gamma^{\Psi^{\rm CAS},\sigma}_{uv} = \langle \Psi^{\rm CAS} | \hat{a}^\dagger_{u\sigma}\hat{a}_{v\sigma}|\Psi^{\rm CAS}\rangle, 
\end{equation}
with all the other blocks (i.e., internal-active and all blocks with at least one external index) vanishing.
The AO spin densities are then given by
\begin{equation}
    \label{eq:CASAOSpinRDM}
    P^{\Psi^{\rm CAS},\sigma} = C\gamma^{{\rm CAS},\sigma} C^T
\end{equation}
and it holds that 
\begin{equation}
    \label{eq:CASProjIneq}
    \Pi^I \leq P^{\Psi^{\rm CAS},\sigma} \leq \Pi^I + \Pi^A.
\end{equation}

Direct optimization methods for ROHF and CASSCF can be divided into two groups, depending on the degrees of freedom used for performing the optimization. 
In the MO formalism, the main variable is the coefficient matrix $C$. As mentioned previously, it is an orthonormal matrix, which can be seen as a point of the orthogonal group $O({\mathcal N}_b)$. In the DM formulation, the main variable is the pair of orthogonal projectors $(\Pi^{\rm I},\Pi^{\rm A})$, which can be identified with a point of the set
\begin{align}\label{def:M_DM}
  \cM_{\rm DM} := \biggl\{(\Pi^{\rm I},\Pi^{\rm A}) &\in \RR^{\cN_{\rm b}\times \cN_{\rm b}}_{\rm sym}\times\RR^{\cN_{\rm b}\times \cN_{\rm b}}_{\rm sym}\mbox{ s.t. }
  (\Pi^{\rm I})^2 = \Pi^{\rm I}, \; (\Pi^{\rm A})^2 = \Pi^{\rm A} \nonumber\\
  &\; {\rm Tr}(\Pi^{\rm I})=N_I, \; {\rm Tr}(\Pi^{\rm A})=N_A, \mbox{ and } \Pi^{\rm I}\Pi^{\rm A} = 0 \biggr\}.
\end{align}
We will see later that the above set has a nice geometrical structure: it can be canonically identified with the flag manifold $\cM_{\rm Flag}:={\rm Flag}(N_{\rm I},N_{\rm I}+N_{\rm A};\R^{\cN_{\rm b}})$.

The passage between MO and DM parameterization is done by the map
\begin{equation}\label{eq:MO->DM}
    \zeta : O({\mathcal N}_b) \ni C \mapsto \left( \Pi^{\rm I}, \Pi^{\rm A}\right) \in \cM_{\rm MO} \qquad \mbox{with } \Pi^{\rm I}, \Pi^{\rm A} \mbox{ given by Eq. \ref{eq:projectors}}.
\end{equation}

The dimension of the MO manifold $O({\mathcal N}_b)$ is $\frac{{\mathcal N}_b({\mathcal N}_b-1)}2$ (i.e. the number of degrees of freedom in an orthogonal matrix), while the dimension of the DM manifold $\cM_{\rm DM}$ can be shown to be $N_I N_A + N_IN_E+N_AN_E$. The discrepancy comes from the fact that rotations that mix orbitals of the same class do not affect the energy. In mathematical terms, this can be formulated as follows: the DM manifold $\cM_{\rm DM}$ can be identified with the quotient of the MO manifold $O({\mathcal N}_b)$ by the group $O(N_I) \times O(N_A) \times O(N_E)$. This identification has very practical consequences on the design of direct optimization algorithms, as will be seen in Section~\ref{sec:optimization_flag}.

The geometric structure described above corresponds to a well known fact in quantum chemistry. In direct optimization implementations, changes in the orbitals are parameterized via a rotation matrix 
\begin{equation}
    \label{eq:UExp}
    U = e^\kappa,
\end{equation}
where $\kappa \in \R^{\cN_{\rm b}\times \cN_{\rm b}}$ is a skew-symmetric matrix with the following block structure:
\begin{equation}
    \label{eq:kappa}
    \kappa = 
    \begin{pmatrix}
        0 & \kappa_{IA} & \kappa_{IE} \\
        - \kappa^T_{IA} & 0 & \kappa_{AE} \\
        - \kappa^T_{IE} & - \kappa^T_{AE} & 0
    \end{pmatrix}
\end{equation} 
The vanishing diagonal blocks, that would mix orbitals belonging to the same class, are the practical translation of the quotient process mentioned above. We further note that the map $\kappa \mapsto Ce^\kappa$, with $\kappa$ as in Eq. \ref{eq:kappa} provides a non-redundant local parametrization of the quotient manifold $O({\mathcal N}_b)/(O(N_I) \times O(N_A) \times O(N_E))$.

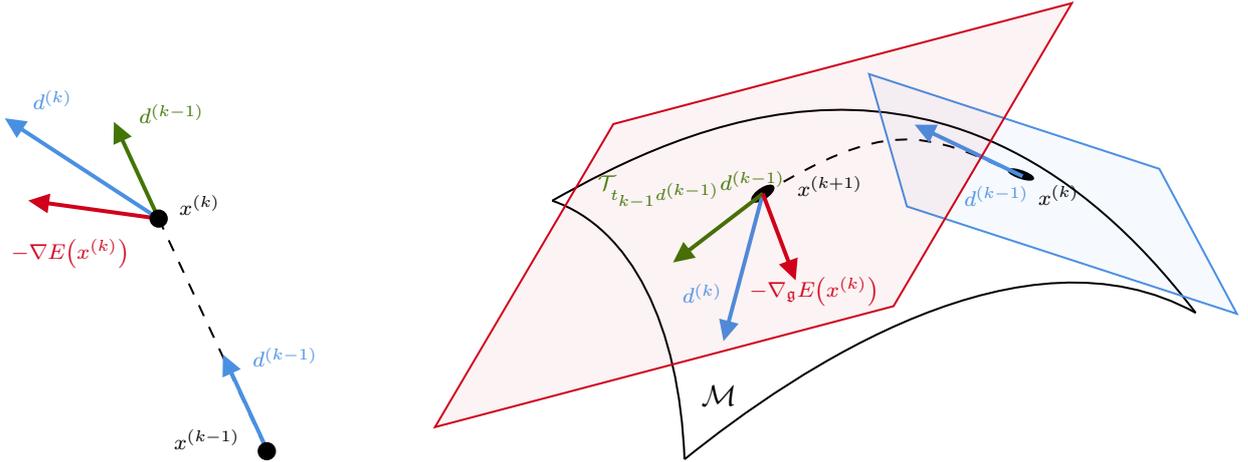
\begin{figure}
    \begin{tabular}{cc}
     \tikzset{every picture/.style={line width=0.75pt}} 

\begin{tikzpicture}[x=0.75pt,y=0.75pt,yscale=-0.9,xscale=0.9]

\draw  [dash pattern={on 4.5pt off 4.5pt}]  (94,88) -- (154,218) ;
\draw [color={rgb, 255:red, 74; green, 144; blue, 226 }  ,draw opacity=1 ][fill={rgb, 255:red, 74; green, 144; blue, 226 }  ,fill opacity=1 ][line width=1.5]    (154.47,218) -- (131.33,167.57) ;
\draw [shift={(129.67,163.93)}, rotate = 65.36] [fill={rgb, 255:red, 74; green, 144; blue, 226 }  ,fill opacity=1 ][line width=0.08]  [draw opacity=0] (11.61,-5.58) -- (0,0) -- (11.61,5.58) -- cycle    ;
\draw [color={rgb, 255:red, 74; green, 144; blue, 226 }  ,draw opacity=1 ][fill={rgb, 255:red, 74; green, 144; blue, 226 }  ,fill opacity=1 ][line width=1.5]    (94,88) -- (12.01,34.13) ;
\draw [shift={(8.67,31.93)}, rotate = 33.31] [fill={rgb, 255:red, 74; green, 144; blue, 226 }  ,fill opacity=1 ][line width=0.08]  [draw opacity=0] (11.61,-5.58) -- (0,0) -- (11.61,5.58) -- cycle    ;
\draw [color={rgb, 255:red, 65; green, 117; blue, 5 }  ,draw opacity=1 ][fill={rgb, 255:red, 74; green, 144; blue, 226 }  ,fill opacity=1 ][line width=1.5]    (94,88) -- (70.87,37.57) ;
\draw [shift={(69.2,33.93)}, rotate = 65.36] [fill={rgb, 255:red, 65; green, 117; blue, 5 }  ,fill opacity=1 ][line width=0.08]  [draw opacity=0] (11.61,-5.58) -- (0,0) -- (11.61,5.58) -- cycle    ;
\draw [color={rgb, 255:red, 208; green, 2; blue, 27 }  ,draw opacity=1 ][fill={rgb, 255:red, 74; green, 144; blue, 226 }  ,fill opacity=1 ][line width=1.5]    (94,88) -- (25.63,78.48) ;
\draw [shift={(21.67,77.93)}, rotate = 7.92] [fill={rgb, 255:red, 208; green, 2; blue, 27 }  ,fill opacity=1 ][line width=0.08]  [draw opacity=0] (11.61,-5.58) -- (0,0) -- (11.61,5.58) -- cycle    ;
\draw  [fill={rgb, 255:red, 0; green, 0; blue, 0 }  ,fill opacity=1 ] (149.47,218) .. controls (149.47,215.5) and (151.5,213.47) .. (154,213.47) .. controls (156.5,213.47) and (158.53,215.5) .. (158.53,218) .. controls (158.53,220.5) and (156.5,222.53) .. (154,222.53) .. controls (151.5,222.53) and (149.47,220.5) .. (149.47,218) -- cycle ;
\draw  [fill={rgb, 255:red, 0; green, 0; blue, 0 }  ,fill opacity=1 ] (89.47,88) .. controls (89.47,85.5) and (91.5,83.47) .. (94,83.47) .. controls (96.5,83.47) and (98.53,85.5) .. (98.53,88) .. controls (98.53,90.5) and (96.5,92.53) .. (94,92.53) .. controls (91.5,92.53) and (89.47,90.5) .. (89.47,88) -- cycle ;

\draw (101,204.4) node [anchor=north west][inner sep=0.75pt]    {$\scriptstyle x^{( k-1)}$};
\draw (10.65,97.72) node [anchor=north west][inner sep=0.75pt]  [color={rgb, 255:red, 208; green, 2; blue, 27 }  ,opacity=1 ] 
{$\scriptstyle -\nabla E\left( x^{( k)}\right)$};
\draw (144.63,158.4) node [anchor=north west][inner sep=0.75pt]  [color={rgb, 255:red, 74; green, 144; blue, 226 }  ,opacity=1 ]  
{$\scriptstyle d^{( k-1)}$};
\draw (81.63,22.4) node [anchor=north west][inner sep=0.75pt]  [color={rgb, 255:red, 65; green, 117; blue, 5 }  ,opacity=1 ] 
{$\scriptstyle d^{( k-1)}$};
\draw (22.63,14.4) node [anchor=north west][inner sep=0.75pt]  [color={rgb, 255:red, 74; green, 144; blue, 226 }  ,opacity=1 ]  
{$\scriptstyle d^{( k)}$};
\draw (104,73.4) node [anchor=north west][inner sep=0.75pt]    {$\scriptstyle x^{( k)}$};

\end{tikzpicture} & \hspace{0.7cm}  \tikzset{every picture/.style={line width=0.75pt}} 

\begin{tikzpicture}[x=0.75pt,y=0.75pt,yscale=-0.9,xscale=0.9]

\draw    (74.75,119.01) .. controls (146.22,78.63) and (296.77,3.22) .. (431.65,181.71) ;
\draw    (148.12,263.55) .. controls (211.97,211.65) and (336.79,128.91) .. (431.65,181.71) ;
\draw    (74.75,119.01) .. controls (118.58,134.15) and (145.26,188.61) .. (148.12,263.55) ;

\draw  [dash pattern={on 4.5pt off 4.5pt}]  (191.55,115.08) .. controls (258.04,75.72) and (278.75,76.55) .. (337.67,107.09) ;
\draw  [fill={rgb, 255:red, 0; green, 0; blue, 0 }  ,fill opacity=1 ] (185.77,117.47) .. controls (186.49,115.36) and (189.65,112.58) .. (192.84,111.26) .. controls (196.04,109.94) and (198.04,110.58) .. (197.33,112.69) .. controls (196.62,114.81) and (193.45,117.59) .. (190.26,118.91) .. controls (187.07,120.23) and (185.06,119.59) .. (185.77,117.47) -- cycle ;
\draw  [fill={rgb, 255:red, 0; green, 0; blue, 0 }  ,fill opacity=1 ] (328.6,103.33) .. controls (326.95,101.81) and (328.34,101.05) .. (331.7,101.61) .. controls (335.05,102.18) and (339.11,103.87) .. (340.76,105.38) .. controls (342.42,106.89) and (341.03,107.66) .. (337.67,107.09) .. controls (334.32,106.53) and (330.26,104.84) .. (328.6,103.33) -- cycle ;
\draw  [color={rgb, 255:red, 74; green, 144; blue, 226 }  ,draw opacity=1 ][fill={rgb, 255:red, 74; green, 144; blue, 226 }  ,fill opacity=0.05 ] (250.47,48.22) -- (411.47,101.22) -- (454.47,182.22) -- (271.47,122.22) -- cycle ;
\draw [color={rgb, 255:red, 74; green, 144; blue, 226 }  ,draw opacity=1 ][fill={rgb, 255:red, 74; green, 144; blue, 226 }  ,fill opacity=1 ][line width=1.5]    (335.6,105.33) -- (279.35,78.28) ;
\draw [shift={(275.75,76.55)}, rotate = 25.68] [fill={rgb, 255:red, 74; green, 144; blue, 226 }  ,fill opacity=1 ][line width=0.08]  [draw opacity=0] (11.61,-5.58) -- (0,0) -- (11.61,5.58) -- cycle    ;
\draw [color={rgb, 255:red, 74; green, 144; blue, 226 }  ,draw opacity=1 ][fill={rgb, 255:red, 74; green, 144; blue, 226 }  ,fill opacity=1 ][line width=1.5]    (191.55,115.08) -- (170.77,193.68) ;
\draw [shift={(169.75,197.55)}, rotate = 284.81] [fill={rgb, 255:red, 74; green, 144; blue, 226 }  ,fill opacity=1 ][line width=0.08]  [draw opacity=0] (11.61,-5.58) -- (0,0) -- (11.61,5.58) -- cycle    ;
\draw [color={rgb, 255:red, 65; green, 117; blue, 5 }  ,draw opacity=1 ][fill={rgb, 255:red, 74; green, 144; blue, 226 }  ,fill opacity=1 ][line width=1.5]    (191.55,115.08) -- (144.92,151.1) ;
\draw [shift={(141.75,153.55)}, rotate = 322.32] [fill={rgb, 255:red, 65; green, 117; blue, 5 }  ,fill opacity=1 ][line width=0.08]  [draw opacity=0] (11.61,-5.58) -- (0,0) -- (11.61,5.58) -- cycle    ;
\draw [color={rgb, 255:red, 208; green, 2; blue, 27 }  ,draw opacity=1 ][fill={rgb, 255:red, 74; green, 144; blue, 226 }  ,fill opacity=1 ][line width=1.5]    (191.55,115.08) -- (208.34,159.81) ;
\draw [shift={(209.75,163.55)}, rotate = 249.42] [fill={rgb, 255:red, 208; green, 2; blue, 27 }  ,fill opacity=1 ][line width=0.08]  [draw opacity=0] (11.61,-5.58) -- (0,0) -- (11.61,5.58) -- cycle    ;
\draw  [color={rgb, 255:red, 208; green, 2; blue, 27 }  ,draw opacity=1 ][fill={rgb, 255:red, 208; green, 2; blue, 27 }  ,fill opacity=0.05 ] (108.69,76.16) -- (362.96,8.6) -- (264.02,177.99) -- (9.75,245.55) -- cycle ;

\draw (342.76,108.78) node [anchor=north west][inner sep=0.75pt]    {$\scriptstyle x^{( k)}$};
\draw (155.82,220.39) node [anchor=north west][inner sep=0.75pt]    {$\mathcal{M}$};
\draw (209.1,103.6) node [anchor=north west][inner sep=0.75pt]    {$\scriptstyle  x^{( k+1)}$};
\draw (301.63,109.4) node [anchor=north west][inner sep=0.75pt]  [color={rgb, 255:red, 74; green, 144; blue, 226 }  ,opacity=1 ]  
{$\scriptstyle d^{( k-1)}$};
\draw (145.47,163.62) node [anchor=north west][inner sep=0.75pt]  [color={rgb, 255:red, 74; green, 144; blue, 226 }  ,opacity=1 ]  
{$\scriptstyle d^{( k)}$};
\draw (182.65,159.72) node [anchor=north west][inner sep=0.75pt]  [color={rgb, 255:red, 208; green, 2; blue, 27 }  ,opacity=1 ]  
{$\scriptstyle -\nabla _{\mathfrak{g}} E\left( x^{( k)}\right)$};
\draw (98.51,100.55) node [anchor=north west][inner sep=0.75pt]  [color={rgb, 255:red, 65; green, 117; blue, 5 }  ,opacity=1 ] 
{$\scriptstyle  \mathcal{T}_{t_{k-1} d^{(k-1)}} d^{( k-1)}$};

\end{tikzpicture}
\end{tabular}
    \caption{Conjugate gradient algorithm in $\R^n$ (left). Riemannian conjugate gradient  (RCG) algorithm on a Riemannian manifold $\mathcal M$ (right).}
    \label{fig:CG_flat_man}
\end{figure}

In standard quantum chemistry direct optimization implementations, a sequence of MO coefficients $\{C^{(k)}\}_{k=0}^{N_{it}}$ is generated starting from an initial guess $C^{(0)}$ such that the sequence of energies $\{E^{(k)}\}_{k=0}^{N_{it}}$ is non increasing and hopefully converges to the ground state energy. The passage from $C^{(k)}$ to $C^{(k+1)}$ is obtained by 
\begin{equation}
    \label{eq:Ck+1}
    C^{(k+1)} = C^{(k)} e^{\kappa^{(k)}}
\end{equation}
where $\kappa^{(k)}$ is the result of some optimization procedure (e.g., steepest descent, or Levenberg-Marquardt second order optimization).  
Eq. \ref{eq:Ck+1} amounts to changing the center of the local parameterization of the quotient manifold. For the steepest descent and Newton optimization methods, the calculation of $\kappa^{(k)}$ relies only on information relative to the point $C^{(k)}$, i.e., it makes no use of the history. Employing an optimization method that does, e.g., non-linear conjugate gradient (CG) or quasi-Newton methods, comes with a complication. 
Let us consider non-linear CG as an example. In the (flat) vector space $\R^n$, the CG descent direction at a given iterate is computed by linearly combining the gradient at the current iterate with the descent direct at the previous direction, see Fig.~\ref{fig:CG_flat_man} (left) and Eq.~\ref{eq:CGA_flat}. On a Riemannian manifold, the gradient and descent direction at a given iterate belong to the tangent space to the manifold at this particular iterate, which changes from iteration to iteration. Therefore, it is not possible to linearly combine tangent vectors at different points, as required by CG, in an obvious way.
The operation of correctly transferring a vector quantity from the tangent space at a given point of the manifold to the tangent space at another point of the manifold is called \emph{transport}, see Fig~\ref{fig:CG_flat_man} (right) and Eq.~\ref{eq:CGA}.

For RHF and RKS, this problem is not apparent because the optimization takes place in a Grassmann manifold, for which the parallel transport map is trivial in the right parameterization (see Eq.~\ref{eq:transport_Grassmann}). This is not the case for the flag manifolds on which ROHF and CASSCF optimization problems are set (see Eq.~\ref{eq:transport_ROHF}). 

\medskip

Given $(\Pi^{\rm I},\Pi^{\rm A}) \in  \cM_{\rm DM}$, there exists a unique (up to an irrelevant global phase) normalized ROHF wavefunction $\Phi_{(\Pi^{\rm I},\Pi^{\rm A})}^{\rm ROHF}$ with maximal spin polarization $S=S_z=\frac{N_A}2$ associated with $(\Pi^{\rm I},\Pi^{\rm A})$: it is the Slater determinant whose spin-1-RDM in the AO basis $(\chi_\mu)$ is given by
Eq. \ref{eq:ROHFSpinRDM}.
The ROHF energy functional 
\begin{equation}\label{eq:ROHF_energy_functional}
E^{\rm ROHF}(\Pi^{\rm I},\Pi^{\rm A}):= \langle \Phi_{(\Pi^{\rm I},\Pi^{\rm A})}^{\rm ROHF}|\widehat H_N|\Phi_{(\Pi^{\rm I},\Pi^{\rm A})}^{\rm ROHF}\rangle
\end{equation}
is therefore a well-defined function of $(\Pi^{\rm I},\Pi^{\rm A})$, and in fact a quadratic function in $\Pi^{\rm I}$ and $\Pi^{\rm A}$. From a geometrical point of view, the ROHF problem is therefore a smooth optimization problem on a flag manifold, for which the energy is quadratic in the density-matrix formalism.

\medskip

CASSCF can be also seen as an optimization problem on a flag manifold. In the spin-collinear approximation, the corresponding CASSCF energy functional can be written as
$$
E^{\rm CAS}(\Pi^{\rm I},\Pi^{\rm A}) = \min_{\Psi \in \cW_{\Pi^{\rm I},\Pi^{\rm A}}^{\rm CAS}}  \langle \Psi|\widehat H_N|\Psi\rangle,
$$
with 
$$
\cW_{\Pi^{\rm I},\Pi^{\rm A}}^{\rm CAS} := \left\{ \Psi  \; \mbox{s.t.} \; \| \Psi \|=1, \; \Pi^{\rm I} \le P^{\Psi,\sigma} \le \Pi^{\rm I}+\Pi^{\rm A}, \; \tr(P^{\Psi,\sigma})=N_I+m_\sigma, \; \sigma=\alpha,\beta \right\}.
$$
Recall that for $A,B \in \R^{\cN_{\rm b} \times \cN_{\rm b}}_{\rm sym}$, $A \le B$ means that $X^TAX \le X^T B X$ for all $X \in \R^{\cN_{\rm b}}$.

\medskip

A very appealing feature of quotient manifolds is that if closed form expressions for parallel transport and geodesics on $\cM$ are known, then closed form expressions for parallel transports on $\cM/G$ can be derived from the ones on $\cM$. The manifold $O(\cN_{\rm b})$ is in fact a Lie group. For this reason, closed form expressions for parallel transport and geodesics on $O(\cN_{\rm b})$ can be constructed from the exponential map.

\section{Optimization on Riemannian manifolds\label{sec:optimization_riemmannian}}

Riemannian optimization (i.e. optimization on manifold endowed with a Riemannian metric) is a major field of computational mathematics with many applications in various areas of science and technology. Several Riemannian optimization libraries are available, in which the most common Riemannian optimization methods are implemented. One of the advantages of using an optimization library is that this allows one to test and compare many different optimization methods with limited development effort. For optimization in the flat space $\R^n$, the user of an optimization library is just asked to provide the code returning the value of the function and its gradient at an input point $x \in \R^n$ (and possibly also, for some methods, a preconditioner and/or the Hessian at $x$). For optimization on a Riemannian manifold $\cM$, the user is asked to provide four pieces of codes returning respectively:
\begin{enumerate}
    \item the value of the function and its Riemannian gradient at an input point $x \in \cM$, (and possibly also, for some methods, a preconditioner and/or the Riemannian Hessian at~$x$);
    \item the value of $\mathfrak g_x(p_x,q_x) \in \R$, where $\mathfrak g_x$ is the Riemannian metric, for an input point $x \in \cM$ and two tangent vectors $p_x,q_x \in T_x\cM$ at point $x$;
    \item the value of $\cR_x(p_x) \in \cM$, where $\cR$ is the chosen retraction, for an input point $x \in \cM$ and a tangent vector $p_x \in T_x\cM$ at point $x$;
    \item the value of $\cT_{p_x}q_x \in T_{\cR_x(p_x)}\cM$, where $\cT$ is the chosen transport, for an input point $x \in \cM$ and two tangent vectors $p_x,q_x \in T_x\cM$ at point~$x$.
\end{enumerate}
 Let us first recall the role of $\mathfrak g$, $\cR$ and $\cT$ in Riemannian optimization algorithms, and illustrate these concepts on the simple example of optimization on the orthogonal group $O({\mathcal N}_b)$. It is well-known that the tangent space to $O({\mathcal N}_b)$ at some $C \in O({\mathcal N}_b)$ is given by
 $$
 T_CO({\mathcal N}_b) = \{ CA, \, A \in \R^{{\mathcal N}_b \times {\mathcal N}_b}_{\rm antisym} \},
 $$
 where $A \in \R^{{\mathcal N}_b \times {\mathcal N}_b}_{\rm antisym}$ is the vector space of ${\mathcal N}_b \times {\mathcal N}_b$ real antisymmetric matrix. The Frobenius inner product 
 $$
 \langle M,N \rangle_{\rm F} :=\tr(M^TN) = \sum_{\mu,\nu=1}^{{\mathcal N}_b} M_{\mu\nu} N_{\mu\nu}
 $$
 on $\R^{{\mathcal N}_b \times {\mathcal N}_b}$ induces a Riemannian metric on $O({\mathcal N}_b)$ defined by
 $$
 \mathfrak g_C(CA,CA')=\tr(A^T A')=-\tr(AA') \quad \mbox{for all } C \in O({\mathcal N}_b) \mbox{ and } CA,CA' \in T_CO({\mathcal N}_b).
 $$
 Consider a smooth function $E: \R^{{\mathcal N}_b \times {\mathcal N}_b} \to \R$. The gradient of $E$ at some point $C \in \R^{{\mathcal N}_b \times {\mathcal N}_b}$ is the unique matrix $\nabla E(C) \in \R^{{\mathcal N}_b \times {\mathcal N}_b}$ such that
 $$
 E(C+\delta C) = E(C) + \langle \nabla E(C),\delta C \rangle_{\rm F} + o(\delta C).
 $$
If $C \in O({\mathcal N}_b)$, the Riemannian gradient of $E$ at $C$ is the matrix $\nabla_{\mathfrak g}E(C) \in T_CO({\mathcal N}_b)$ obtained by orthogonally projecting $\nabla E(C)$ on $T_CO({\mathcal N}_b)$ for the Frobenius inner product. Its expression is given by 
$$
\nabla_{\mathfrak g}E(C) = \frac 12 C \left( C^T \nabla E(C) - \nabla E(C)^T C \right).
$$

\begin{figure}
        \centering
        \begin{tikzpicture}[x=0.75pt,y=0.75pt,yscale=-1,xscale=1]

\draw    (113.45,88.05) .. controls (188.45,49.5) and (346.45,-22.5) .. (488,147.91) ;
\draw    (190.45,226.05) .. controls (257.45,176.5) and (388.45,97.5) .. (488,147.91) ;
\draw    (113.45,88.05) .. controls (159.45,102.5) and (187.45,154.5) .. (190.45,226.05) ;

\draw  [color={rgb, 255:red, 74; green, 144; blue, 226 }  ,draw opacity=0.72 ][fill={rgb, 255:red, 74; green, 144; blue, 226 }  ,fill opacity=0.09 ] (325.45,3.05) -- (495,67.91) -- (253.45,144.05) -- (146.45,41.05) -- cycle ;
\draw  [dash pattern={on 4.5pt off 4.5pt}]  (172.45,139.5) .. controls (227.45,69.5) and (253.45,44.5) .. (357,53.91) ;
\draw  [fill={rgb, 255:red, 0; green, 0; blue, 0 }  ,fill opacity=1 ] (179.52,129.17) .. controls (179.45,126.71) and (181.66,122.63) .. (184.45,120.05) .. controls (187.24,117.47) and (189.56,117.37) .. (189.63,119.83) .. controls (189.7,122.28) and (187.49,126.36) .. (184.7,128.94) .. controls (181.91,131.52) and (179.59,131.62) .. (179.52,129.17) -- cycle ;
\draw [color={rgb, 255:red, 208; green, 2; blue, 27 }  ,draw opacity=1 ][line width=1.5]    (188.45,85.5) .. controls (176.89,97.21) and (177.71,104.26) .. (183.15,119.47) ;
\draw [shift={(184.45,123.05)}, rotate = 249.75] [fill={rgb, 255:red, 208; green, 2; blue, 27 }  ,fill opacity=1 ][line width=0.08]  [draw opacity=0] (11.61,-5.58) -- (0,0) -- (11.61,5.58) -- cycle    ;
\draw [color={rgb, 255:red, 74; green, 144; blue, 226 }  ,draw opacity=1 ][fill={rgb, 255:red, 74; green, 144; blue, 226 }  ,fill opacity=1 ][line width=1.5]    (274,56.71) -- (195.82,80.73) ;
\draw [shift={(192,81.91)}, rotate = 342.92] [fill={rgb, 255:red, 74; green, 144; blue, 226 }  ,fill opacity=1 ][line width=0.08]  [draw opacity=0] (11.61,-5.58) -- (0,0) -- (11.61,5.58) -- cycle    ;
\draw [color={rgb, 255:red, 74; green, 144; blue, 226 }  ,draw opacity=1 ][fill={rgb, 255:red, 208; green, 2; blue, 27 }  ,fill opacity=1 ][line width=1.5]    (277,55.71) -- (330.79,40.11) ;
\draw [shift={(334.63,39)}, rotate = 163.83] [fill={rgb, 255:red, 74; green, 144; blue, 226 }  ,fill opacity=1 ][line width=0.08]  [draw opacity=0] (11.61,-5.58) -- (0,0) -- (11.61,5.58) -- cycle    ;
\draw  [fill={rgb, 255:red, 0; green, 0; blue, 0 }  ,fill opacity=1 ] (265.65,56.71) .. controls (267.61,54.91) and (273.39,53.46) .. (278.55,53.46) .. controls (283.72,53.46) and (286.32,54.91) .. (284.35,56.71) .. controls (282.39,58.5) and (276.61,59.96) .. (271.45,59.96) .. controls (266.28,59.96) and (263.68,58.5) .. (265.65,56.71) -- cycle ;

\draw (280.55,56.86) node [anchor=north west][inner sep=0.75pt]    {$x^{( k)}$};
\draw (199.63,53.4) node [anchor=north west][inner sep=0.75pt]  [color={rgb, 255:red, 74; green, 144; blue, 226 }  ,opacity=1 ]  {$t_{k} d^{( k)}$};
\draw (407.63,57.11) node [anchor=north west][inner sep=0.75pt]  [color={rgb, 255:red, 74; green, 144; blue, 226 }  ,opacity=1 ]  {$T_{x}\mathcal{M}$};
\draw (197.63,182.11) node [anchor=north west][inner sep=0.75pt]    {$\mathcal{M}$};
\draw  [color={rgb, 255:red, 0; green, 0; blue, 0 }  ,draw opacity=1 ]  (318,166) -- (475,166) -- (475,199) -- (318,199) -- cycle  ;
\draw (320,170.4) node [anchor=north west][inner sep=0.75pt]    {$x^{( k+1)} =\mathcal{R}_{x^{( k)}}\left( t_{k} d^{( k)}\right)$};
\draw (303,15.4) node [anchor=north west][inner sep=0.75pt]  [color={rgb, 255:red, 74; green, 144; blue, 226 }  ,opacity=1 ]  {$-\nabla _{\mathfrak{g}} E$};
\draw (139,79.4) node [anchor=north west][inner sep=0.75pt]  [color={rgb, 255:red, 208; green, 2; blue, 27 }  ,opacity=1 ]  {$\mathcal{R}_{x^{( k)}}$};
\draw (181.52,132.57) node [anchor=north west][inner sep=0.75pt]    {$x^{( k+1)}$};

\end{tikzpicture}
        \caption{Riemannian steepest descent (RSD) algorithm.}
        \label{fig:SD_man}
    \end{figure}
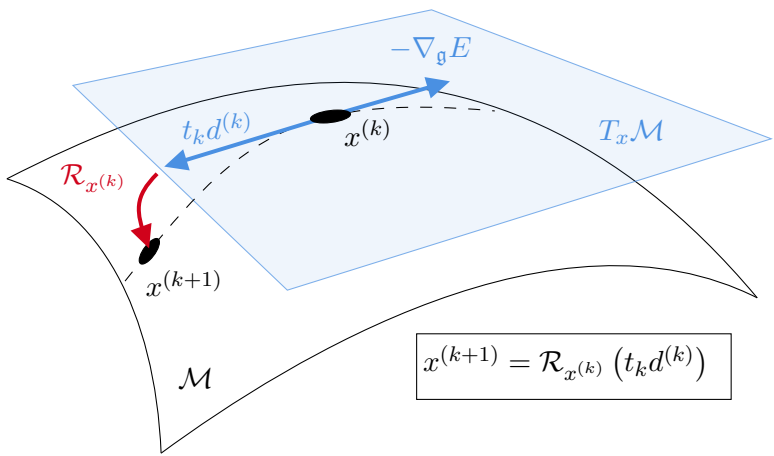

A retraction $\cR$ on a manifold $\cM$ is a map $\cR :T\cM \to \cM$ such that for all $x \in \cM$ the restriction $\cR_x:T_x\cM \to \cM$ of $\cR$ to $T_x\cM$ satisfies for $p_x \in T_x\cM$
    \begin{equation}\label{eq:retraction}
    \cR_x(p_x)=x + p_x + o(p_x).
    \end{equation}
    Among other things (see below), retractions are used to map straight lines, or more generally paths, drawn on the vector space $T_x\cM$ onto paths drawn on the curved manifold $\cM$. As a matter of example, the iterates of the fixed-step gradient (also called steepest descent) algorithm are defined by
    \begin{equation}\label{eq:SDA}
    d^{(k)} = - \nabla_{\mathfrak g} E(x^{(k)}), \quad x^{(k+1)}=\cR_{x^{(k)}}(t d^{(k)}), \qquad \mbox{(fixed-step steepest descent)}
    \end{equation}
    for a chosen fixed step $t > 0$.
    In words, starting from a point $x^{(k)} \in \cM$, the descent direction $d^{(k)}$ is chosen equal to the opposite of the gradient, which is the steepest descent direction for infinitesimal length steps, a step $t d^{(k)}$ is made in this direction, and finally, the vector $t d^{(k)} \in T_{x^{(k)}}\cM$ is mapped back to a point of the manifold thanks to the retraction (see Fig.~\ref{fig:SD_man}). 

    It follows from  \eqref{eq:retraction} that in the limit of small step lengths, we have
    $$
    x^{(k+1)}=\cR_{x^{(k)}}(t d^{(k)}) = x^{(k)} + t d^{(k)} + o(t)
    $$
    where the remainder term $o(t)$ can be interpreted as a correction due to the curvature of the manifold $\cM$.

Among all possible retractions on a Riemannian manifold $\cM$, one is canonical: it is the one defined from the geodesics, called the exponential map, and denoted by ${\rm Exp}$. In the  case of $O({\mathcal N}_b)$, the exponential map has a simple closed expression and is related to the usual exponential of matrices:
\[
{\rm Exp}_C(CA) = Ce^A, \quad A \in \R^{{\mathcal N}_b \times {\mathcal N}_b}_{\rm antisym} \quad \mbox{(exponential map on $O({\mathcal N}_b)$).}
\]

In Riemannian optimization, vector transports are used in particular to combine the descent directions of previous iterates. Let us further elaborate on this point, that was qualitatively discussed in Section \ref{sec:theory}. In the standard conjugate gradient algorithm in $\R^n$, the descent direction $d^{(k)}$ at iterate $x^{(k)}$ is a linear combination of $-\nabla E(x^{(k)})$ and $d^{(k-1)}$, the previous descent direction:
\begin{equation} \left| \begin{array}{l}
    d^{(k)} = - g^{(k)}+ \beta_k d^{(k-1)},  \\
      x^{(k+1)} = x^{(k)} + t_k d^{(k)}, \end{array}\right. \quad \mbox{(nonlinear CG algorithm in $\R^n$)}  \label{eq:CGA_flat} 
    \end{equation}
with $g^{(k)} = \nabla E(x^{(k)})$ and either
\begin{align*}
    & \beta_{k}^{\rm FR}:= \frac{\Vert g^{(k)} \Vert^2}{\Vert g^{(k-1)}\Vert^2} & & \mbox{(Fletcher-Reeves)}, \\
     & \beta_{k}^{\rm PR}:= \frac{{g^{(k)}}^T(g^{(k)}-g^{(k-1)})}{\Vert g^{(k-1)}\Vert^2} & & \mbox{(Polak-Ribi\`ere)},\\
    &  \beta_{k}^{\rm HS} := \frac{{g^{(k)}}^T(g^{(k)}-g^{(k-1)})}{(g^{(k)}-g^{(k-1)})^T d^{(k)}} & & \mbox{(Hestenes-Stiefel)}.
\end{align*}
The step length $t_k \in \R$ is obtained by a line search technique such as Armijo, Wolfe, or Hager-Zhang \cite{hager2006algorithm} linesearch.

This idea cannot be directly used in optimization on manifolds, because $-\nabla E(x^{(k)})$ and $d^{(k-1)}$ belong to different vector spaces, namely $T_{x^{(k)}}\cM$ and $T_{x^{(k-1)}}\cM$ respectively. Before being combined with $-\nabla E(x^{(k)})$ to form the new descent direction $d^{(k)}$, the vector $d^{(k-1)}$ must be transported from $T_{x^{(k-1)}}\cM$ to $T_{x^{(k)}}\cM$, using a transport map $\cT$. A transport map takes as input two vectors $p_x,q_x$ of the tangent space $T_x\cM$ at point $x$, and returns a vector $\cT_{p_x}q_x$ of the tangent space at point $\cR_x(p_x)$ (compatibility condition with the retraction $\cR$). The map $(p_x,q_x) \mapsto \cT_{p_x}q_x$ is linear in the variable $q_x$, and satisfies the consistency relation
    $\cT_0 q_x=q_x$. We thus have
    \begin{equation} \left| \begin{array}{l}
    d^{(k)} = - g^{(k)}+ \beta_k \cT_{t_{k-1}d^{(k-1)}}d^{(k-1)},  \\
      x^{(k+1)} = \cR_{x^{(k)}}(t_k d^{(k)}), \end{array}\right. \quad \mbox{(Riemannian CG algorithm)}  \label{eq:CGA} 
    \end{equation}
with either 
\begin{align*}
    & \beta_{k}^{\rm RFR}:= \frac{{\mathfrak g}_{x^{(k)}}(g^{(k)},g^{(k)})}{\mathfrak g_{x^{(k-1)}}(g^{(k-1)},g^{(k-1)})} & & \mbox{(Riemannian Fletcher-Reeves)}, \\
     & \beta_{k}^{\rm RPR}:= \frac{\mathfrak g_{x^{(k)}}(  g^{(k)},g^{(k)} - \cT_{t_{k-1}d^{(k-1)}}(g^{(k-1)}) ) }{\mathfrak g_{x^{(k)}}(g^{(k-1)},g^{(k-1)})} & & \mbox{(Riemannian Polak-Ribi\`ere)},\\
     & \beta_{k}^{\rm RHS} := \frac{\mathfrak g_{x^{(k)}}(g^{(k)},g^{(k)} - \cT_{t_{k-1}d^{(k-1)}}(g^{(k-1)}))}{\mathfrak g_{x^{(k)}}(g^{(k)},\cT_{t_{k-1}d^{(k-1)}}(d^{(k-1)})) - g_{x^{(k-1)}}(g^{(k-1)},d^{(k-1)})} & & \mbox{(Riemannian Hestenes-Stiefel)}.
\end{align*}
Together with $x^{(k)} = \cR_{x^{(k-1)}}(t_{k-1} d^{(k-1)})$, the fact that $\cT_{p_x}q_x \in T_{\cR_x(p_x)}\cM$ (comptabilitly with the retraction) ensures that $\cT_{t_{k-1}d_{(k-1)}}d^{(k-1)}$ belongs to $T_{x^{(k)}}\cM$.

Among all transports compatible with the exponential map ${\rm Exp}$ associated with the metric $\mathfrak g$, one is canonical: it is the parallel transport associated with the Levi-Civita connection of the metric $\mathfrak g$. 
For the example of $O({\mathcal N}_b)$, this parallel transport has an extremely simple form
$$
\cT_{CA}(CB) = Ce^A B \quad \mbox{(parallel transport on $O({\mathcal N}_b)$)}.
$$

\section{Optimization on Grassmann and flag manifolds\label{sec:optimization_flag}}

In RHF and RKS models, the state of the system is described by a point of the Grassmann manifold
$$
{\rm Gr}(N,{\mathcal N}_b)  \cong \underbrace{\left\{ P \in \R^{{\mathcal N}_b \times {\mathcal N}_b}_{\rm sym} \mbox{ s.t. } P^2=P, \; \tr(P)=N \right\}}_{\rm DM \; formalism} \cong \underbrace{O({\mathcal N}_b)/(O(N) \times O({\mathcal N}_b-N))}_{\rm MO \; formalism}.
$$
In the DM formalism, the Grassmann manifold is parameterized by the matrix $P$ of the orthogonal projector on the vector space spanned by the doubly-occupied MO. In the MO formalism, it is represented by an orthogonal matrix $C \in O({\mathcal N}_b)$, the first $N$ columns of $C$ corresponding to the $N$ doubly-occupied orbitals, and the last ${\mathcal N}_b-N$ ones to the virtual orbitals. The gauge invariance in the MO formulation is taken into account by quotienting $O({\mathcal N}_b)$ by the group $O(N) \times O({\mathcal N}_b-N)$ (occupied-occupied and virtual-virtual rotations, repsectiely).

\medskip

Likewise, in the ROHF model and the outer CASSCF minimization problem, the state is represented by a point of the flag manifold  
$$
{\rm Flag}(N_I,N_I+N_A,{\mathcal N}_b)  \cong \cM_{\rm DM}  \cong \underbrace{O({\mathcal N}_b)/(O(N_I) \times O(N_A) \times O(N_E))}_{\rm MO \; formalism}.
$$

In both cases, the MO formalism involves the quotient of the orthogonal group $O({\mathcal N}_b)$ by a closed subgroup ($O(N) \times O({\mathcal N}_b-N)$ for RHF/RKS, $O(N) \times O(N_I) \times O(N_A) \times O(N_E)$ for ROHF/CASSCF). As a consequence, the closed form expressions for the canonical retraction and parallel transport on $O({\mathcal N}_b)$ can be translated into closed form expressions for canonical retraction and parallel transport on the quotient manifold \cite{Edelman1998,Ye2022,absil2008optimization,boumal2023}, leading to the following formulae:
\begin{itemize}
    \item RHF/RKS setting: 
    \begin{align} 
&    \mbox{tangent space at $C$} \quad \cong  \left\{ C\kappa ,   \quad   \kappa= \begin{pmatrix}
        0 & \kappa_{OV}  \\
        - \kappa^T_{OV} & 0 
    \end{pmatrix}\right\}, \label{eq:TM_Grassmann} \\
    & \mbox{metric:} \qquad 
    \mathfrak g_{C}(C\kappa,C\kappa') = \tr( \kappa^T \kappa') = 2 \tr( \kappa_{OV}^T\kappa_{OV}'), \label{eq:matrix_Grassmann} \\
    & \mbox{exponential map (canonical retraction):} \qquad 
    \cR_{C}(C\kappa) = C e^{\kappa}, \label{eq:retraction_Grassmann} \\
    & \mbox{parallel transport:} \qquad 
    \cT_{C\kappa}(C\kappa') = Ce^\kappa \kappa', \label{eq:transport_Grassmann}
    \end{align} 
    \item ROHF/CASSCF setting:
        \begin{align} 
&    \mbox{tangent space at $C$} \quad \cong  \left\{ C\kappa ,   \quad   \kappa= \begin{pmatrix}
        0 & \kappa_{IA} & \kappa_{IE}  \\
        - \kappa_{IA}^T & 0 &  \kappa_{AE}  \\ 
        - \kappa_{IE}^T & - \kappa_{AE}^T & 0
    \end{pmatrix} \right\}, \label{eq:TM_ROHF} \\
    & \mbox{metric:} \qquad 
    \mathfrak g_{C}(C\kappa,C\kappa') = \tr( \kappa^T \kappa') 
    \label{eq:matrix_ROHF} \\
    & \mbox{exponential map (canonical retraction):} \qquad 
    \cR_{C}(C\kappa) = C e^{\kappa}, \label{eq:retraction_ROHF} \\
    & \mbox{parallel transport:} \qquad 
    \cT_{C\kappa}(C\kappa') = Ce^\kappa e^{-\phi_\kappa}(\kappa'), \label{eq:transport_ROHF}
\end{align}
where $\phi_\kappa : \mathfrak K \to \mathfrak K$ is the linear operator on the vector space
$$
\mathfrak K= \left\{ \kappa= \begin{pmatrix}
        0 & \kappa_{IA} & \kappa_{IE}  \\
        - \kappa_{IA}^T & 0 &  \kappa_{AE}  \\ 
        - \kappa_{IE}^T & - \kappa_{AE}^T & 0
    \end{pmatrix} \right\} 
$$
defined by
\begin{align*}
\phi_\kappa (\kappa') &= \frac 12 {\rm Proj}_{\mathfrak K}\left( [ \kappa,\kappa'] \right) \\
&= \frac 12 \begin{pmatrix}
        0 & - \kappa_{IE} [\kappa_{AE}']^T + \kappa_{IE}' \kappa_{AE}^T & \kappa_{IA}\kappa_{AE}'- \kappa_{IA}'\kappa_{AE}  \\
        \kappa_{AE}' \kappa_{IE}^T - \kappa_{AE} [\kappa_{IE}']^T & 0 &  - \kappa_{IA}^T \kappa_{IE}' + [\kappa_{IA}']^T \kappa_{IE} \\ 
   - [\kappa_{AE}']^T  \kappa_{IA}^T+ \kappa_{AE}^T [\kappa_{IA}']^T   & [\kappa_{IE}']^T \kappa_{IA}  - \kappa_{IE}^T \kappa_{IA}'  & 0
    \end{pmatrix},
\end{align*}
and 
\begin{equation}\label{eq:transport_exp_operator}
    e^{-\phi_\kappa} = \sum_{n=0}^{+\infty} \frac{(-1)^n}{n!} \underbrace{(\phi_\kappa \circ \cdots \circ \phi_\kappa)}_{n \; \rm times}.
\end{equation}

\end{itemize}

In computational codes, it is convenient to represent a tangent vector $C\kappa$ by the block $\kappa_{OV} \in \R^{N \times ({\mathcal N}_b-N)}$ for RHS/RKS, and the blocks $(\kappa_{IA},\kappa_{IE},\kappa_{AE}) \in \R^{N_I \times N_A} \times \R^{N_I \times N_E} \R^{N_A \times N_E}$ for ROHF/CASSCF. It follows from Eq. \ref{eq:transport_Grassmann}, that in this representation the parallel  transport for RHF/RKS is the identity operator. This is not the case for ROHF/CASSCF where the transport of the block matrix $\kappa'$ is done by the map $e^{-\phi_\kappa}(\kappa')$, which transforms and mixes the IA/IE/AE blocks of $\kappa'$. Let us note however that in the special case when the transported vector $C \kappa'$ is collinear to the vector $C\kappa$ along which it is transported, then the transport formula Eq. \ref{eq:transport_ROHF} dramatically simplifies. Indeed, we then have $[\kappa,\kappa']=0$, and therefore $e^{-\phi_\kappa}(\kappa') = \kappa'$. This occurs for the Riemannian conjugate gradient method (see Eq. \ref{eq:CGA}), but not for quasi-Newton methods such as BFGS.

\section{Numerical Results}\label{sec:numerical_results}

In this section, we analyze the performance of Riemannian optimization algorithms for solving the ROHF and CASSCF minimization problems for a few selected test cases. Let us first provide some implementation details.

\paragraph{General implementation.} We focus specifically on Riemannian steepest descent (RSD), nonlinear conjugate gradient (RCG) and low-memory Broyden-Fletcher-Goldfarb-Shanno (R-LBFGS) methods, all endowed with preconditioning. We refer to Refs \cite{absil2008optimization,boumal2023} for  general introductions to Riemannian optimization methods.  
Our code is structured as follows: first, we implemented the RSD, RCG, and R-LBFGS optimization routines in the MO formalism within a Julia \cite{bezanson2017julia} package which is then interfaced with PySCF \cite{sun2020recent} for ROHF and CFOUR \cite{matthews2020coupled} for CASSCF calculations. These software handle the operations specific
to the ROHF and CASSCF models, including the generation of AO basis sets and initial guess MOs, 
the computation of electronic integrals, and the evaluation of energies and Frobenius gradients.

In our Julia package, we use for all methods the exponential retraction \eqref{eq:retraction_ROHF} and parallel transport \eqref{eq:transport_ROHF} as outlined in the previous section. For parallel transport, the exponential operator
is computed by truncating the series \eqref{eq:transport_exp_operator} so that the Frobenius norm of the last term falls below numerical precision.

Our implementation of RCG is based on {\it Algorithm 1} in Boumal et al. \cite{boumal2015low} with Polak-Ribière coefficient $\beta^{\rm RPR}$ as above. For R-LBFGS, we implemented {\it Algorithm~2} in Huang et al. \cite{huang2015broyden}. For CASSCF, we use the inverse diagonal of the Hessian as preconditioner. In the case of ROHF, we tested two different preconditioners. The first one is the modified inverse diagonal Hessian as discussed in Ref \cite{nottoli_black-box_2021}. The other is detailed in Appendix. Our results are presented for the second choice of preconditioner that showcased the best performance for our test case. 

All methods use Hager-Zhang \cite{hager2006algorithm} linesearch as implemented in the LineSearches.jl \cite{mogensen2018optim} Julia package. Computations are considered to have reached convergence when
the Frobenius norm of the Riemannian gradient reaches $10^{-5}$. 
Comprehensive implementation details are available in our publicly accessible GitHub repository.\footnote{\url{https://github.com/LaurentVidal95/ROHFToolkit}}

\medskip

\paragraph{Details specific to the R-LBFGS implementation.} At each iteration, the R-LBFGS method constructs an approximation $B$ of the inverse Hessian using a certain number $m$ of vectors stored in memory from previous iterations, through an iterative procedure  \cite{huang2015broyden}. In addition to preconditioning, it is important to note that the performance of R-LBFGS is influenced by the selection of the maximum depth $m_{\rm max}$, the initial guess $B_0$ for the approximate inverse Hessian in the iterative process and the choice of restart strategies, which determines the iterations at which the history is reset. For both ROHF and CASSCF, we define $B_0=\gamma\,{\rm Id}$ with $\gamma$ as in Ref~\cite{huang2015broyden}. If $P$ is the preconditioner, the preconditioned version of R-LBFGS is obtained by replacing $B_0$ with $\hat{\gamma} P$, with $\hat{\gamma}$ as in Ref~\cite{de2018nonlinearly}. 

We experimented two restart strategies which depend on the preconditioning. For the first one, called “dynamic” R-LBFGS, the diagonal Hessian used for preconditioning is updated at each iteration. The history is reset whenever the direction obtained from the R-LBFGS quasi-newton system is not a descent direction. For the second method called “fixed” R-LBFGS, we use the same preconditioner $P={\rm diag}({\rm Hess}^{(0)})^{-1}$, corresponding to the inverse diagonal Hessian for the guess orbitals, at each iteration until the inverse diagonal Hessian at current point, ${\rm diag}({\rm Hess}^{(k)})^{-1}$, deviates too much from $P$. When this happens, the history is reset, and the procedure is reinitialized with $P={\rm diag}({\rm Hess}^{(k)})^{-1}$. When using the preconditioner described in appendix for ROHF, we applied the dynamic strategy.

\subsection{ROHF}

\begin{figure}[h!]
    \centering
    \begin{tabular}{cc}
         \includegraphics[width=0.488\textwidth]{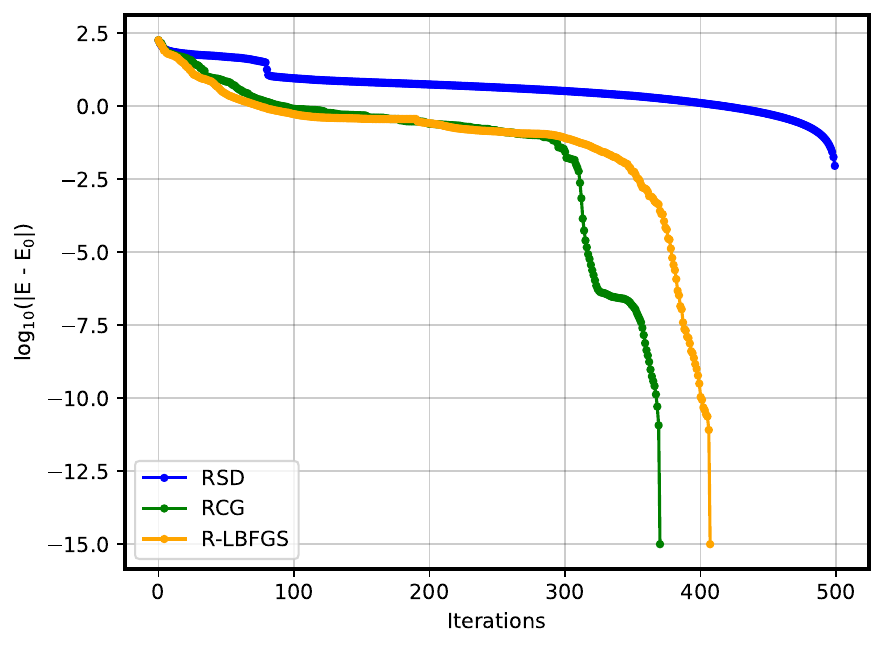}
         &\includegraphics[width=0.47\textwidth]{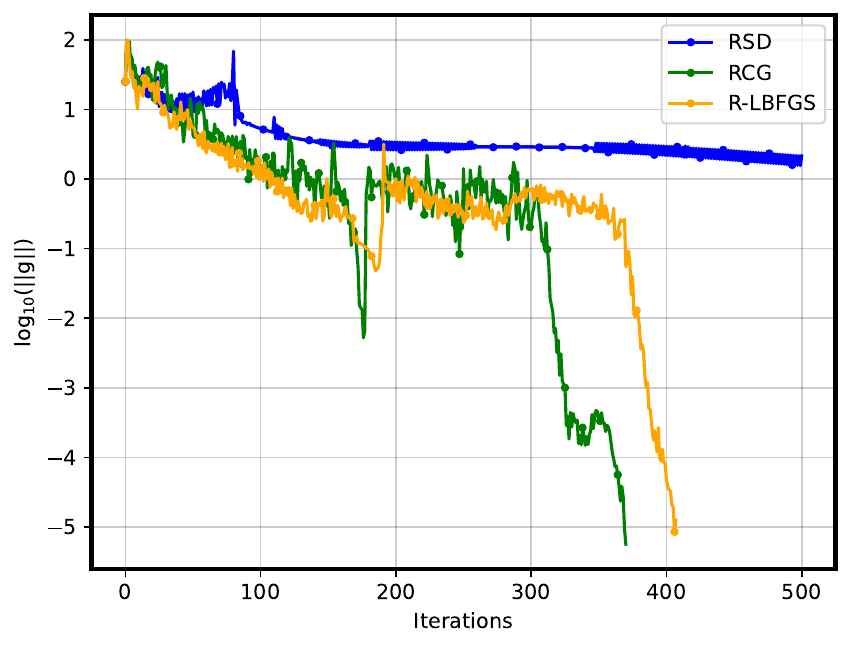}
    \end{tabular}
    \caption{Calculations from core initial guess. On the left, energy difference with respect to the converged energy along the iterations. On the right, Frobenius norm of the Riemannian gradient along the iterations. Only the first 500 iterations of RSD are shown on the graph for readability.}
    \label{fig:conv_ROHF_bad}
\end{figure}

\begin{figure}[h!]
    \centering
    \begin{tabular}{cc}
         \includegraphics[width=0.477\textwidth]{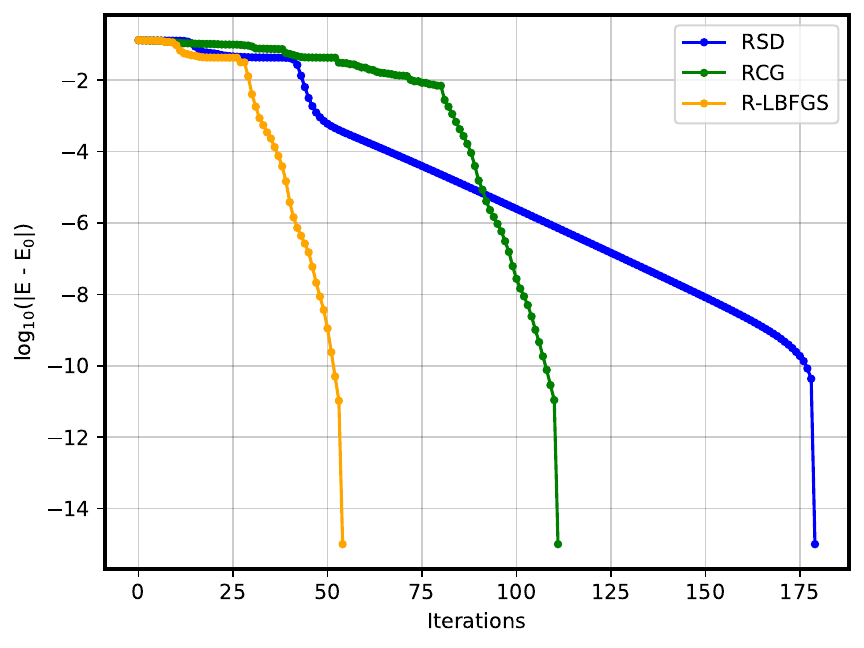}
         &\includegraphics[width=0.47\textwidth]{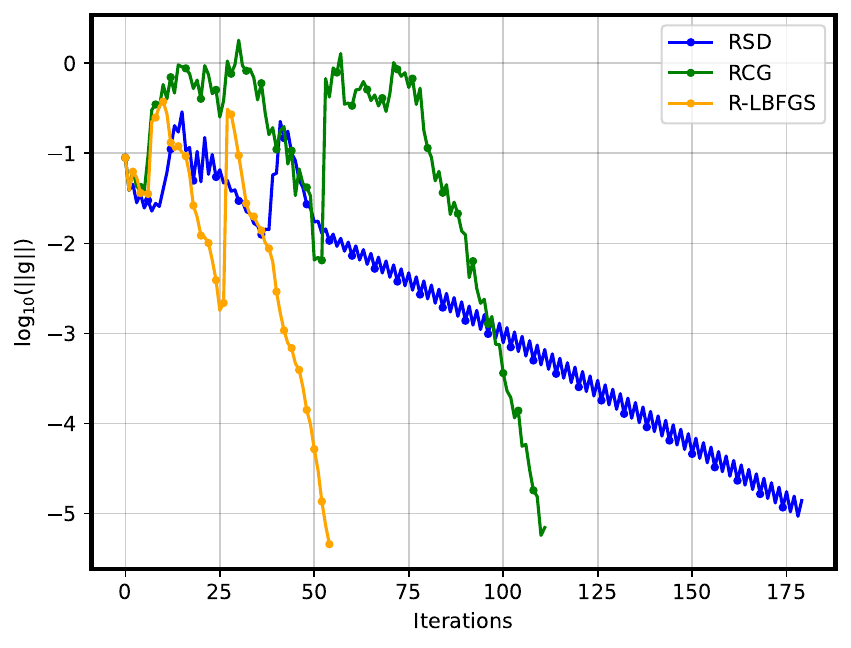}
    \end{tabular}
    \caption{Calculations from a better initial guess, 1 $E_h$ from the expected energy. On the left, energy difference with respect to the converged energy along the iterations. On the right, Frobenius norm of the Riemannian gradient along the iterations.}
    \label{fig:conv_ROHF_good}
\end{figure}

For ROHF we tested the three aforementioned Riemannian optimization methods on Ti$_2$O$_4$ in its $D_{\rm 2h}$ geometry, using Dunning's \textit{cc}-pVTZ basis set \cite{Dunning1989,Kendall1992}. This system is employed as a template for addressing SCF convergence issues\footnote{See \url{https://www.scm.com/doc/ADF/Examples/SCF_Ti2O4.html}} in the Amsterdam Density Functional (ADF) quantum-chemistry package \cite{te2001chemistry}. In order to compare the performance of Riemannian algorithms in different convergence regimes, calculations were started from both a core initial guess (Fig.~\ref{fig:conv_ROHF_bad}) and a guess closer to a minimum (Fig.~\ref{fig:conv_ROHF_good}). The second guess is obtained by a standard SCF+DIIS method for ROHF, with Guest and Saunders coefficients \cite{plakhutin2014canonical}, stopped when the Frobenius norm of the gradient reaches $10^{-1}$.

\begin{figure}[h!]
    \centering
    \includegraphics[width=0.488\textwidth]{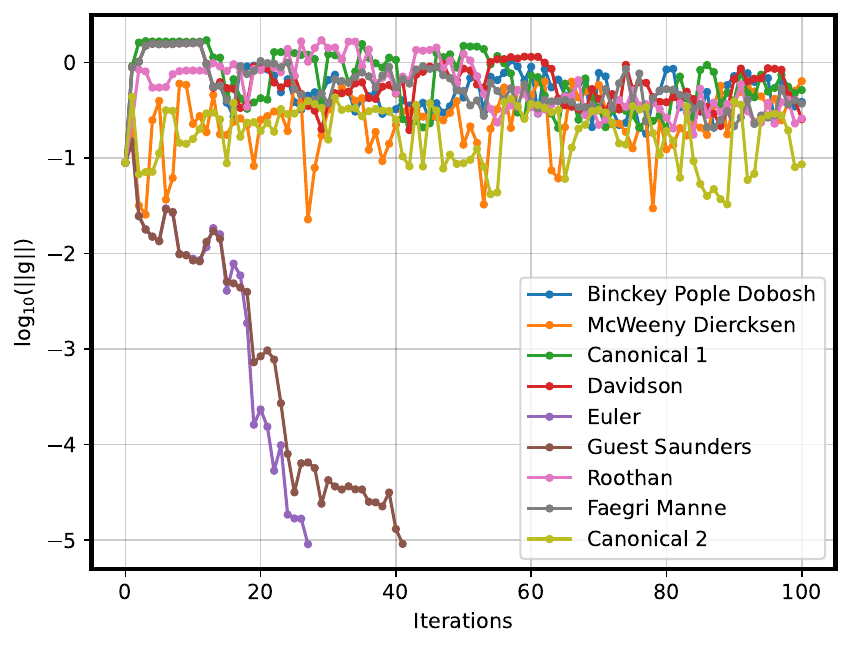}
    \caption{Comparison of standard SCF+DIIS methods for ROHF as listed in Ref \cite{plakhutin2014canonical}, for
    Ti$_2$O$_4$ in \textit{cc}-pVTZ basis set, starting from our good initial guess
    about 1 Hartree away from the expected energy. 
    The vertical axis shows the Frobenius norm of the Riemannian gradient along the iterations.
    Only two methods yield convergence.}
    \label{fig:SCFs_ROHF_good}
\end{figure}

In both cases, all three methods provide stable convergence toward a local minimum of the energy. In the optimal scenario, a finely tuned SCF+DIIS method outperforms the Riemannian optimization methods we have tested. However, the performance and stability of SCF routines for ROHF are notoriously sensitive to the choice of method and acceleration parameters, as illustrated in Fig.~\ref{fig:SCFs_ROHF_good}. 
On the other hand, the direct minimization methods described in this paper have the advantage of offering robust convergence, which is a valuable feature in terms of user's time and effort. 

We were unfortunately not able to make a direct comparison to the GDM algorithm \cite{Dunietz2002} due to our lack of access to the code or to the fine details of the implementation. Nevertheless, a simple-minded test performed with the free trial version of Q-Chem \cite{shao2015advances} using default parameters showed that GDM and GDM+DIIS exhibit similar performances to our RCG implementation on Ti$_2$O$_4$.

An important point that needs to be discussed here is what minimum the various algorithms converged to. 
Our proof-of-concept implementation does not enforce point-group symmetry, so our calculations were performed in the $C_{\rm 1}$ group. 

Our calculations converged to two different minima, one at $-1996.191285~E_h$, which was systematically obtained when starting from the core guess, and one at $-1996.179398~E_h$, obtained when starting from the better guess.
We also looked at the lowest triplets for each Irrep enforcing symmetry using a quadratically convergent ROHF implementation \cite{nottoli_black-box_2021}, and we determined that the lowest triplet is the $B_{\rm 2u}$ state at $-1996.142005$~$E_h$. The stability analysis of such solution reveals however that a lower energy, symmetry-broken solution exists. We therefore conclude that the two solutions found with our Riemannian optimization algorithms are two lower-energy symmetry-broken solutions. Whether a lower-energy, symmetry broken-solution is desirable or not depends on the aims of the study, and ultimately on the user; however, we note that our algorithms can be generalized to enforce point group symmetry, which we plan to do in the future.
\color{black}

\subsection{CASSCF}
The CASSCF method has been tested by running calculations on a subgroup of the benchmark set used in Ref.~\cite{Menezes2016} and Ref.~\cite{Nottoli2021} using Pople's 6-31G* basis set \cite{hehre1972self}. Convergence properties of direct minimization algorithms were compared against two well-established CASSCF optimization algorithms namely Super CI (SCI)~\cite{roos1980_sci, siegbahn1981_sci,malmqvist1990,kollmar2019_ptsci,angeli2002_ptsci} and the norm-extended optimization (NEO)~\cite{Jensen1983,Jensen1986}, the latter one being a genuine second-order algorithm. 
All computations were carried out using two different choices, to simulate, respectively, a troublesome scenario where the calculation starts relatively far away from the converged result, and an ideal starting point that should be close to the final minimum. For the former scenario, we use canonical restricted Hartree-Fock (RHF) orbitals, while as a good starting point we exploit unrestricted natural orbitals (UNO)~\cite{Pulay1988,Toth2020}

\begin{figure}
    \centering
    \includegraphics[width=0.93\textwidth]{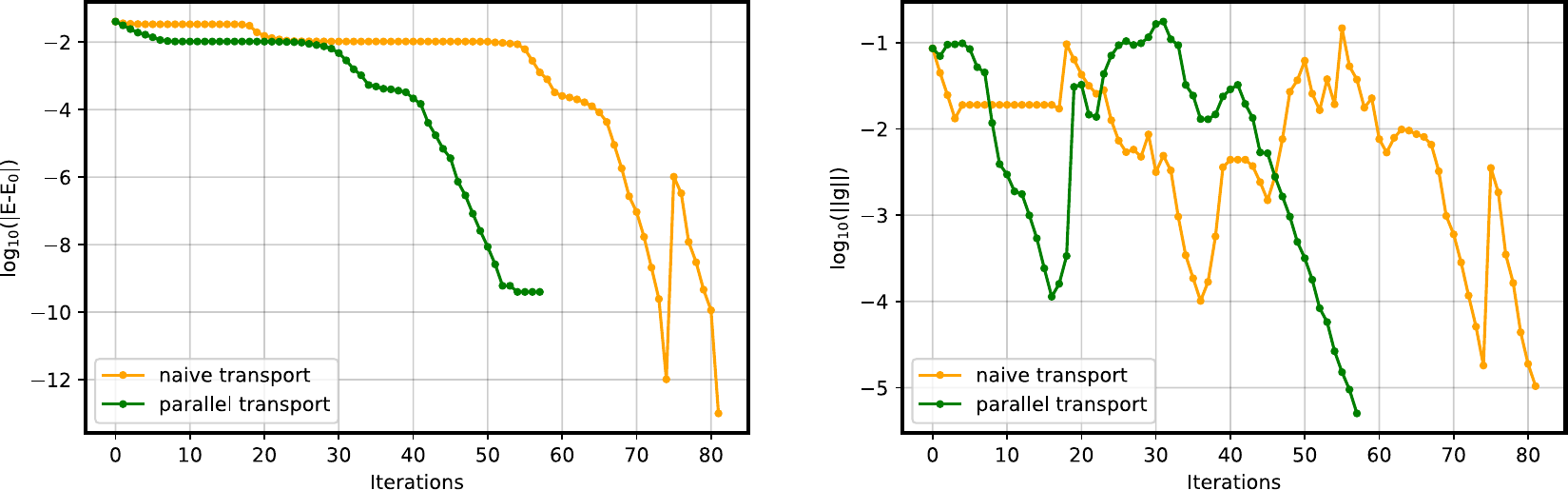}
    \caption{On the left, energy difference with respect to the converged energy along the iterations for a proper R-LBFGS implementation that makes use of the parallel transport (green curve) and a more basic one (orange curve). On the right, Frobenius norm of the gradient along the iterations. }
    \label{fig:conv_trend}
\end{figure}

We report in Tab.~\ref{tab:avg_it} for each algorithm the average number of iterations required to converge CASSCF. The number of iterations for each tested system is reported in the supporting information. As expected, the values related to the RHF guess are systematically higher than the ones related to the UNO guess. Moreover, we notice that the numbers related to the RHF guess show a high variability as indicated by the large standard deviation, thus being strongly system dependent. The average number of iterations for all direct minimization methods with the exception of RSD is comparable with and in some cases outperforms the ones of SCI.
\begin{table}[h!]
    \centering
    \begin{tabular}{lcccc}
         \toprule
         Algorithm & $\langle \text{It.}\rangle^{\text{RHF}}$ & $\sigma^{\text{RHF}}$ & $\langle \text{It.}\rangle^{\text{UNO}}$ & $\sigma^{\text{UNO}}$\\
         \midrule
         RSD              &  115.9   &  112.9   &      32.3   &     12.5 \\
         RCG              &   31.9   &   11.6   &      15.4   &      3.0 \\
         R-LBFGS(dyn)      &   37.9   &   11.3   &      19.2   &      3.5 \\
         R-LBFGS(fix)      &   35.6   &   15.2   &      19.1   &      3.7 \\
         SCI(DIIS)        &   38.4   &   26.0   &      14.9   &      6.2 \\
         SCI(no DIIS)     &   61.5   &   25.1   &      21.6   &      9.9 \\
         NEO              &   12.2   &    1.8   &       5.1   &      0.5 \\
         \bottomrule
    \end{tabular}
    \caption{Average number of iterations ($\langle\text{It.}\rangle$) and standard deviation ($\sigma$) for each tested algorithm starting with two different guess orbitals, namely restricted Hartree-Fock (RHF) and unrestricted natural orbitals (UNO).}
    \label{tab:avg_it}
\end{table}
We conclude this section by looking more in detail at one example, namely, pyridine using a standard CAS(6,6) wavefunction. In Fig.~\ref{fig:conv_trend}, we compare the convergence behavior of the R-LBFGS as implemented in the present study with a naive implementation that simply translates the gradient from previous points, without parallel transport. In both cases, we start from canonical orbitals. We note that while both implementations get stuck for a while on a plateau, the R-LBFGS overcomes it in a few iterations and then converges smoothly. On the contrary, the naive LBFGS implementation exhibits worse convergence, with very small and even positive slope steps at the beginning of the optimization. This demonstrates the importance of properly accounting for parallel transport. 

Another important point concerns the actual solution obtained. When starting from a poor guess such as the one given by canonical orbitals, without any kind of manual selection for the active space, the optimizer may get stuck in local minima that are ultimately related by orbitals swapping between the inactive and active domains. Using once again pyridine as a test case, we observe that all the algorithms converged to the same local minimum (-246.766857 $E_h$) with the exception of SCI that converged to a higher minimum (-246.756489 $E_h$). 

\begin{figure}
    \centering
    \includegraphics[width=0.6\textwidth]{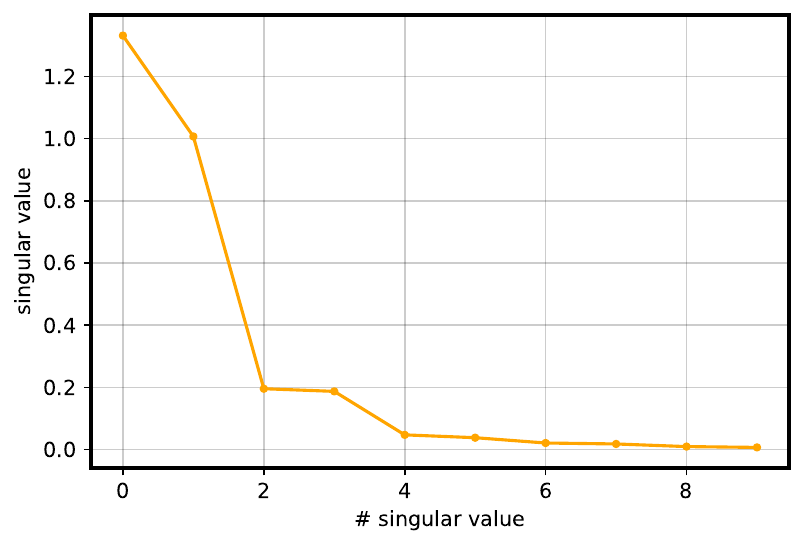}
    \caption{First ten singular values of the difference matrix between the active AO-based one-body density matrices stemming from two calculations that reached different local minima.}
    \label{fig:svd}
\end{figure}

These two minima are characterized by different converged active orbitals. This difference can be assessed simply by visual inspection or by checking the singular values of the difference between the active part of the one-body density matrix in the AO basis for the two calculations. If the converged active orbitals were (almost) the same, we would observe (almost) vanishing singular values. On the contrary, as depicted in Fig.~\ref{fig:svd}, two orbitals are completely different between the two results.

We also note that by manually tuning the DIIS parameters, we managed to achieve convergence to the lower-energy solution using SCI. Again, this underlines the robustness of the Riemannian optimization algorithms, which is valuable in terms of user time and effort.

\section{Conclusions and perspectives}
In this contribution, we have explored the use of Riemannian optimization methods on the flag manifold to optimize restricted-open and complete active space self-consistent field wavefunctions. After discussing the geometry of the problem, we have reviewed the general aspects of Riemannian optimization and its application to the aforementioned chemical problem. We have then compared various algorithms to traditional ones. The Riemannian optimization methods illustrated in this work all show robust convergence properties, and do so without requiring the user to finely tune the parameters that control the optimization. Even in the naive implementation presented here, they demonstrate that they can be competitive with other traditional implementations in terms of number of iterations, and thus overall computational cost. Nevertheless, this is just a proof-of-concept study, for which several further developments are required. First, the overall underwhelming performance of the Riemannian Quasi-Newton L-BFGS method for CASSCF, which is expected to outperform conjugate gradient, as observed for ROHF when starting from a good guess (see Fig. \ref{fig:conv_ROHF_good}), can be explained by the basic preconditioner and initialization of the inverse Hessian we use. Finding better preconditioners and inverse Hessian approximations is not a straightforward task, and requires further attention. We also have not investigated optimal parameters, including exploring different line-search and restarting strategies, which would require to run extensive numerical tests, but that can greatly improve the overall performance of the methods.

In conclusion, we believe that Riemannian optimization is a valuable addition to the SCF optimization toolbox for ROHF and CASSCF, and that further exploration of the use of such techniques is worthy of attention.

\section*{Acknowledgements} This project has received funding from the European Research Council (ERC) under the European Union's Horizon 2020 research and innovation program (grant agreement EMC2 No. 810367). F.L. and T.N. acknowledge financial support from ICSC-Centro Nazionale di Ricerca in High Performance Computing, Big Data, and Quantum Computing, funded by the European Union -- Next Generation EU -- PNRR, Missione 4 Componente 2 Investimento 1.4.

\bibliography{bibliography}
\bibliographystyle{unsrt}

\appendix

\section{A variant of the diagonal Hessian preconditioner for ROHF}

Let $C\in O({\mathcal N}_b)$ and $\kappa\in \mathfrak{K}$ such that $C\kappa\in T_CO({\mathcal N}_b)$.
The ROHF Hessian applied to $C\kappa$ is given by
\begin{equation}
    \mathcal{L}_C(C\kappa) = C
     \begin{pmatrix}
        0 & X & Y  \\
        -X^T & 0 & Z\\
        -Y^T & -Z^T & 0
    \end{pmatrix}
\end{equation}
where the matrices $X\in\mathbb{R}^{N_I\times N_A}$, $Y\in\mathbb{R}^{N_I\times N_E}$ and $Z\in\mathbb{R}^{N_A\times N_E}$ are defined by
\begin{equation}
    \begin{array}{l}
         X = 2(\kappa_{IA}(F_I - F_A)_{AA} - (F_I - F_A)_{II}\kappa_{IA}) +
         \kappa_{IE}(2F_I - F_A)_{EA} + (F_I - 2F_A)_{IE}\kappa_{AE}^T\\
         \hspace{6cm} + (2J(\lambda_1) - K(\lambda_1))_{IA} + J(\lambda_2)_{IA}\\
         Y = \kappa_{IA}(4F_I - 2F_A)_{AE} + 4(\kappa_{IE}(F_I)_{EE} - (F_I)_{II}\kappa_{IE}) - 2((F_I)_{IA} + (F_A)_{IA})\kappa_{AE}\\
         \hspace{6cm} + 4(2J(\lambda_1) - K(\lambda_1))_{IE} + 2(2J(\lambda_2) - K(\lambda_2))_{IE}\\
         Z = \kappa_{IA}^T(2F_I - F_A)_{IE} - 2(F_I+F_A)_{AI}\kappa_{IE} + 4(\kappa_{AE}(F_A)_{EE} - (F_A)_{AA}\kappa_{AE})\\
         \hspace{6cm} + 2(2J(\lambda_1) - K(\lambda_1))_{AE} + 2(2J(\lambda_2) - K(\lambda_2))_{AE}.
    \end{array}
\end{equation}
In the above expression, we adopted the following conventions: the operators $J$ and $K$ are the standard exchange and Coulomb operators, $F_I$ and $F_A$ are the internal and active Fock matrices, defined for all
$\Pi_I=C_I^{ }C_I^T$ and $\Pi_A= C_A^{ }C_A^{T}$ by
\begin{equation*}
\begin{array}{l}
     F_I = \displaystyle h + 2J(\Pi_I) + J(\Pi_A) - K(\Pi_I) - \frac{1}{2}K(\Pi_A)\\[0.15cm]
     F_A = \displaystyle \frac{1}{2}(h + 2J(\Pi_I) + J(\Pi_A) - K(\Pi_I) - K(\Pi_A))
\end{array}
\end{equation*}
and the matrices $\lambda_1,\lambda_2\in\RR^{{\mathcal N}_b\times {\mathcal N}_b}_{\rm sym}$ are given by
\begin{equation}
    \lambda_1 = \begin{pmatrix}
        0 & \kappa_{IA} & \kappa_{IE}\\
        \kappa_{IA}^T & 0 & 0\\
        \kappa_{IE}^T & 0 & 0
    \end{pmatrix}\quad\mbox{and}\quad
    \lambda_2 = \begin{pmatrix}
        0 & -\kappa_{IA} & 0\\
        -\kappa_{IA}^T & 0 & \kappa_{AE}\\
        0 & \kappa_{AE}^T & 0
    \end{pmatrix}.
\end{equation}
Each block $X$, $Y$ and $Z$ can be decomposed as the sum of two
terms, the first one, denoted by $(\widetilde X,\widetilde Y, \widetilde Z)$, being simple to compute in terms of internal
and active Fock operators, and the second one, denoted by $(\Omega_X,\Omega_Y,\Omega_Z)$, being more costly to compute:
\begin{equation}\label{eq:ROHF_hessian_blocs}
    \begin{array}{l}
         X = 2(\kappa_{IA}(F_I - F_A)_{AA} - (F_I - F_A)_{II}\kappa_{IA}) + \Omega_X := \widetilde{X} + \Omega_X\\
         Y = 4(\kappa_{IE}(F_I)_{EE} - (F_I)_{II}\kappa_{IE}) + \Omega_Y := \widetilde{Y} + \Omega_Y\\
         Z = 4(\kappa_{AE}(F_A)_{EE} - (F_A)_{AA}\kappa_{AE}) + \Omega_Z := \widetilde{Z} + \Omega_Z
    \end{array}.
\end{equation}
In our implementation, we define a preconditioned direction
$C\kappa_{\rm prec}$ as a solution to the linear system 
\begin{equation}\label{eq:ROHF_quasi_newton}
    \widetilde{L}_C(C\kappa_{\rm prec}) = C\kappa.
\end{equation}
involving the approximate Hessian
\begin{equation}
    \widetilde{\mathcal{L}}_C(C\kappa) = C \begin{pmatrix}
        0 & \widetilde{X} & \widetilde{Y}\\
        -\widetilde{X}^T & 0 & \widetilde{Z}\\
        -\widetilde{Y}^{T} & -\widetilde{Z}^T & 0
    \end{pmatrix}.
\end{equation}
The advantage of formulation \eqref{eq:ROHF_hessian_blocs} is that the lowest
eigenvalue of $\widetilde L_C$ can be estimated with respect to $F_I$ and $F_A$, which allows to apply a shift when $\widetilde{\mathcal{L}}$ is not positive definite (which is expected when starting far from a minimum).
In addition, the system \eqref{eq:ROHF_quasi_newton}
reads as three Sylvester matrix equations, that can be solved using standard LAPACK optimized routines.

\end{document}





We report in Table~\ref{tab:it_rhf} and \ref{tab:it_uno} the number of iterations required to converge CASSCF for various molecular systems using direct minimazion methods and standard CASSCF optimization algorithms. In Table~\ref{tab:it_rhf} the RHF canonical orbitals were used as initial guess for the CASSCF optimization while in Table~\ref{tab:it_uno} we used unrestricted natural orbitals (UNO).
\begin{table}[h!]
    \centering
    \resizebox{\columnwidth}{!}{%
    \begin{tabular}{lcccccccc}
        \toprule
        system &   active space &  RSD &    RCG     & R-LBFGS(dyn) &   R-LBFGS(fix) &   SCI(DIIS) &   SCI(no DIIS) &   NEO \\
        \midrule
        benzene     & (6,6) & 57    &  57  &    35   &  29   &    17   &   76   &     12 \\
        biphenyl    & (12,12) & 34    &  27  &    40   &  22   &    31   &   53   &     13 \\
        catechol    & (6,6) & 41    &  20  &    24   &  19   &    19   &   48   &      8 \\
        adrenaline  & (6,6) & 32    &  20  &    31   &  24   &    18   &   42   &     12 \\
        pyridine    & (6,6) &236    &  47  &    57   &  56   &    27   &    0   &     15 \\
        azulene     & (10,10) & 31    &  21  &    25   &  24   &    22   &   38   &     12 \\
        fluorene    & (12,12) & 93    &  24  &    29   &  27   &    24   &  109   &     10 \\
        niacina     & (6,6) & 39    &  23  &    32   &  22   &    75   &   60   &     11 \\
        pyridoxal   & (8,8) & 54    &  28  &    35   &  44   &    64   &    0   &     13 \\
        tryptophan  & (8,8) &172    &  39  &    43   &  51   &    44   &    0   &     13 \\
        niacinamide & (6,6) &184    &  40  &    37   &  38   &     0   &    0   &     14 \\
        indole      & (8,8) & 34    &  21  &    32   &  22   &    25   &   36   &     10 \\
        dopamine    & (6,6) & 46    &  25  &    27   &  27   &    27   &   46   &     14 \\
        naphthalene & (10,10) & 89    &  41  &    51   &  44   &    37   &  110   &     12 \\
        pyridoxin   & (6,6) &384    &  32  &    39   &  50   &    38   &    0   &     13 \\
        2Me2HSdiox  & (4,4) & 68    &  25  &    39   &  31   &    28   &   59   &     11 \\
        nicotine    & (6,6) &376    &  53  &    68   &  76   &   119   &    0   &     15 \\
        \bottomrule
    \end{tabular}
    }
    \caption{Number of iterations required to converge the CASSCF wave function for different systems using RHF canonical orbitals as starting guess.}
    \label{tab:it_rhf}
\end{table}

\begin{table}[h!]
    \centering
    \resizebox{\columnwidth}{!}{%
    \begin{tabular}{lcccccccc}
        \toprule
        system &    active space & RSD &    RCG     & R-LBFGS(dyn) &   R-LBFGS(fix) &   SCI(DIIS) &   SCI(no DIIS) &   NEO \\
        \midrule
        benzene     & (6,6) & 23   & 12   &   16   &  15   &    8   &    13  &    5 \\
        biphenyl    & (12,12) & 19   & 12   &   16   &  15   &    8   &    12  &    6 \\
        catechol    & (6,6) & 24   & 14   &   17   &  17   &   10   &    18  &    5 \\
        adrenaline  & (6,6) & 51   & 21   &   26   &  26   &   15   &    26  &    5 \\
        pyridine    & (6,6) & 22   & 12   &   16   &  15   &    9   &    12  &    5 \\
        azulene     & (10,10) & 15   & 11   &   16   &  15   &    5   &     6  &    4 \\
        fluorene    & (12,12) & 20   & 13   &   16   &  16   &    9   &    12  &    6 \\
        niacina     & (6,6) & 56   & 20   &   25   &  25   &   11   &    32  &    6 \\
        pyridoxal   & (8,8) & 39   & 18   &   21   &  20   &   16   &    34  &    5 \\
        tryptophan  & (8,8) & 51   & 20   &   26   &  26   &   23   &    35  &    5 \\
        niacinamide & (6,6) & 24   & 14   &   16   &  16   &   22   &    28  &    5 \\
        indole      & (8,8) & 29   & 15   &   18   &  18   &   20   &    20  &    5 \\
        dopamine    & (6,6) & 32   & 16   &   20   &  21   &   19   &    19  &    5 \\
        naphthalene & (10,10) & 25   & 14   &   17   &  18   &   11   &    11  &    5 \\
        pyridoxin   & (6,6) & 32   & 15   &   19   &  20   &   23   &    25  &    5 \\
        2Me2HSdiox  & (4,4) & 49   & 17   &   21   &  21   &   23   &    42  &    5 \\
        nicotine    & (6,6) & 38   & 17   &   21   &  21   &   22   &    23  &    5 \\
        \bottomrule
    \end{tabular}
    }
    \caption{Number of iterations required to converge the CASSCF wave function for different systems using unrestricted natural orbitals (UNO) as starting guess.}
    \label{tab:it_uno}
\end{table}